\documentclass{tac}
\usepackage[all]{xy}
\usepackage{amstext, amssymb,hyperref,graphicx}
\usepackage{latexcad}
\newcommand{\e}{\mathbf{E}}
\newcommand{\cspge}{Cspn(Graph/|\e|)}
\newcommand{\cspmge}{Cspn(MonGraph/|\e|)}
\newcommand{\cospe}{cospan(\e)}

\usepackage{amsmath}

\begin{document}

\copyrightyear{2007}
\title{ Calculating colimits compositionally}
\author{R. Rosebrugh, N. Sabadini, and R.F.C Walters}
\thanks{The authors gratefully acknowledge financial support from 
the  Universit\'a dell'Insubria and the Italian Government PRIN project ART ({\em Analisi di sistemi di Riduzione mediante sistemi di Transizione}), and the Canadian NSERC.}
\address{Department of Mathematics and Statistics\\
Mt. Allison University, Sackville,\\ New Brunswick, Canada\\[3pt]
and\\
Dipartimento di Scienze delle Cultura,
Politiche e dell'Informazione,\\
Universit\`a dell'Insubria, Italy\\
}

\eaddress{rrosebrugh@mta.ca, nicoletta.sabadini@uninsubria.it,\\ 
robert.walters@uninsubria.it}
\date{\\ Universit\`a degli Studi dell'Insubria, Como, Italy}
\maketitle

\begin{abstract}We show how finite limits and colimits can be calculated compositionally  using the algebras of spans and cospans, and give as an application a proof of the Kleene Theorem on regular languages.
\end{abstract}

\tableofcontents

\section{ Introduction}

\emph{In Computer Science:}

The state spaces of systems are often described by finite limits or colimits in a category $\e$ parametrized by a graph $G$ which describes the underlying geometry of the system.
It is desirable that there is also an algebraic description, so that the limit or colimit is described by an expression, rather than geometrically.

This goes back to the beginnings of computer science, where (i) a program may be described either by a flow chart 
       (goto's), or program text (while) (B\"ohm-Jacopini), (ii) a language may be specified by an automaton or
       an expression  (Kleene). And of course it is present in innumerable areas of computer science (Petri nets versus process algebras, wysiwig versus markup, graph versus term rewriting, etc. etc.) and mathematics.

\medskip

\noindent\emph{In Category Theory:}

Finite limits and colimits are parametrized by graphs; that is, geometrically.
We show that they can also be described by expressions in an algebra.
As an application we prove Kleene's theorem.

\bigskip The algebra in which finite limits and colimits in $\e$ may be expressed compositionally is a appropriate structure on spans and cospans  in the category. This fact is a partial explanation for the algebra  of spans and cospans introduced in \cite{KSW97b},\cite{KSW00b} and developed in various papers, such as \cite{KSW00a},\cite{KSW02},\cite{RSW04}.    

This note is an expanded version of  a lecture to Category Theory 2007, Carvoeiro, Portugal, 18th June 2007. The actual slides of the lecture are available at the \href{http://www.mat.uc.pt/%7Ecateg/ct2007/\#1}{CT2007 web site}. A more detailed version with full proofs is in preparation \cite{RSW08}.

\section{What algebra?}
Assume now \(\mathbf{}\mathbf{E}\) is a category with finite colimits. What is the algebra in which finite colimits in \(\mathbf{E}\) can be described by expressions?

\noindent It is \(cospan(\e)\), considered as a symmetric monoidal category in which each object has a commutative separable algebra structure.
We call a category with such a structure \emph{wscc} (well-supported compact closed \cite{W87}).
 To be precise:

\begin{definition}

A \emph{commutative separable algebra \cite{CW87}} in a symmetric monoidal category is an object \(A \text{ equipped 
with four arrows}\)
$$!:I\longrightarrow A,\:\ \nabla:A\otimes A\longrightarrow A,\:\emph{\text{\textexclamdown}} :A\longrightarrow I,\:\Delta:A\longrightarrow A\otimes A$$
such that
\((A,\triangledown,!) \) forms a commutative monoid, \((A,\vartriangle,\emph{\text{\textexclamdown}})    \)
 forms a cocommutative comonoid, and the following three axioms hold\[(1_A\otimes\nabla)
 (\Delta\otimes1_{A})=\Delta\nabla=(\nabla\otimes1_{A})(1_A\otimes\Delta),\]
\[\nabla\Delta=1_{A}.
\]
\end{definition}

\noindent We can draw a picture of the these three extra axioms, namely:

\noindent\centerline{{\tt\setlength{\unitlength}{0.60pt}
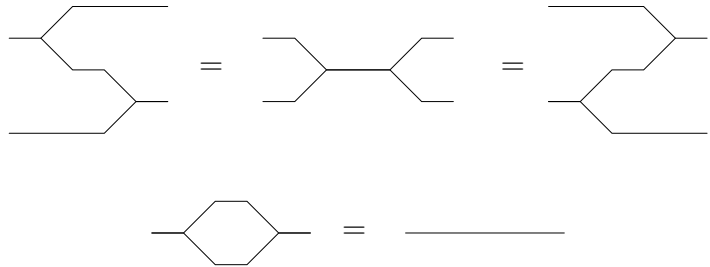
\begin{picture}(440,120)
\thinlines
\put(310,37){$=$}
\put(120,37){$=$}
\path(0,0)(60,0)(80,20)(100,20)
\path(0,60)(20,60)(40,40)(60,40)(80,20)(100,20)
\path(20,60)(40,80)(100,80)
\path(160,20)(180,20)(200,40)(240,40)(260,20)(280,20)
\path(160,60)(180,60)(200,40)(240,40)(260,60)(280,60)
\path(340,20)(360,20)(380,0)(420,0)(440,0)
\path(340,20)(360,20)(380,40)(400,40)(420,60)(440,60)
\path(340,80)(400,80)(420,60)
\end{picture}}}
 
\noindent\centerline{{\tt\setlength{\unitlength}{0.60pt}
\begin{picture}(260,80)
\thinlines
\put(120,17){$=$}
\path(0,20)(20,20)(40,0)(60,0)(80,20)(100,20)
\path(0,20)(20,20)(40,40)(60,40)(80,20)(100,20)
\path(160,20)(260,20)
\end{picture}}}

\medskip
The wscc structure induces a self-dual compact closed structure on the category, and we denote the units and counits of this structure as 

$$\eta_A:I\to A\otimes A\: (=\Delta\ \cdot\ !),\: \varepsilon_A:A\otimes A\to I\:\ (=\text{\textexclamdown}\cdot\nabla). \ $$
For some background to these axioms see also \cite{K04}.

\subsection{The wscc structure on Span and cospan categories}
We will describe the wscc structure on  $\cospe $
for $\e$ a finitely cocomplete category -- the dual structure on $span(\e)$ will then be clear.

An object of $\cospe$ is an object of $\e$; an arrow of $\cospe$ from $A$ to $B$ is an isomorphism class of cospans from $A$ to $B$; that is, of pairs of arrows $$\alpha_1,\alpha_2:A\to R\leftarrow B.$$
\noindent  We will use the notation $\alpha_1,\alpha_2;A\longleftrightarrow B$ to distinguish cospans from arrows in $\e$.  However given any arrow $f:A\to B$ there is a special cospan denoted $\vec f=f,1_B:A\longleftrightarrow B$ and $f^{o}=1_{B},f:B\longleftrightarrow A$. Composition of cospans is by pushout. Now to describe the wscc structure of $\cospe$. The monoidal structure is sum. The special arrows $$!:I\longrightarrow A,\:\ \nabla:A\otimes A\longrightarrow A,\:\text{\textexclamdown}:A\longrightarrow I,\:\Delta:A\longrightarrow A\otimes A$$ are (using $\nabla$ both for the codiagonal in $\e$ and the structure in $\cospe$, and similarly overloading the symbol $!$) $$!=\vec!:0\longleftrightarrow A,\:\ \nabla=\vec \nabla:A+ A\longleftrightarrow A,\:\text{\textexclamdown}=!^o:A\longleftrightarrow 0,\:\Delta=\nabla^o:A\longleftrightarrow  A+A.$$

\subsection{$\cspge$}
Let $Graph$ be the category of finite graphs, let $|\e|$ be the underlying graph (possibly infinite) of $\e$.
Consider $Graph/|\e|$, the category with objects \emph{diagrams} in $\e$, and morphisms \emph{compatible graph morphisms}.  Then  $\cspge$ is the full subcategory of $cospan(Graph/|\e|)$ whose objects are \emph{discrete} diagrams in $\e$. 

\emph{Notice that colimits in this category are calculated as in $Graph$, and are unrelated to colimits in $\e$.}  

We will denote diagrams using set-theoretical notation; for example $\xymatrix{
\{A\ar@<1ex>@{->}[r]^-{f}\ar@<-1ex>@{->}[r]_-{g}&B\}}$ denotes the diagram with two parallel arrows. 
 
It will be useful to introduce a way of picturing arrows in $\cspge $ (engineering notation). Represent the objects in the centre graph of the cospan as points, and arrows in the centre as components with one input (to the left) and one output (to the right) joined to those points which are the domain and codomain of the arrow.  Represent the graph morphisms of the cospan as input and output wires of the whole picture. 
\begin{example}\label{compeqn} 

The following cospan of diagrams $$\xymatrix{\{A\}\ar[r]&
\{A\ar@<1ex>@{->}[r]^-{f}&B\ar@<+1ex>@(ul,ur)^{h}\ar[r]^{k}\ar@<+1ex>@{->}[l]^-{g}&C\}&\{C\}\ar[l]}.$$
This cospan could be pictured as

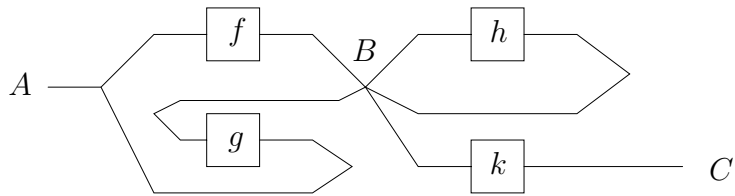
\centerline{{\tt\setlength{\unitlength}{1.00pt}
\begin{picture}(155,80)
\thinlines
\put(-15,37){$A$}
\put(115,50){$B$}
\put(250,5){$C$}
\put(68,18){$g$}
\put(68,58){$f$}
\put(167,57){$h$}
\put(167,7){$k$}
\path(0,40)(20,40)(40,60)(60,60)(60,70)(80,70)(80,50)(60,50)(60,60)
\path(120,40)(110,35)(50,35)(40,30)(50,20)(60,20)(60,30)(80,30)(80,10)(60,10)(60,20)
\path(80,20)(100,20)(115,10)(100,0)(40,0)(20,40)
\path(120,40)(140,60)(160,60)(160,70)(180,70)(180,50)(160,50)(160,60)
\path(80,60)(100,60)(120,40)
\path(180,60)(200,60)(220,45)(200,30)(140,30)(120,40)
\path(120,40)(140,10)(160,10)(160,20)(180,20)(180,0)(160,0)(160,10)
\path(180,10)(240,10)
\end{picture}}}
\end{example}

\section{The Theorem}

\noindent Taking colimit of diagrams in $\e$ induces a functor
 $$colim:Graph/|\e|\to \e.$$
\noindent 

 \begin{theorem}\label{theorem}The functor $colim:Graph/|\e|\to \e$ extends to a functor
      $$ colim : \cspge \to  \cospe$$
  which preserves the wscc structure. 
 \end{theorem}
         
\noindent The definition of the extended colimit is just applying colimit to cospans. It is straightforward that this colim preserves the constants of wscc structure of $\cspge$, and that colim of the cospan of the diagram  $\xymatrix{
\{A\}\ar@{->}[r]^-{}&\{A\ar[r]^{f}&B\}&\{B\}\ar[l]}$ is $\vec f:A\to B$.
 The fact that colim preserves the tensor is also clear. What remains to prove is the fact that colim is a functor -- we outline the proof below. 

Another special case of colim is worth remarking. Consider a cospan in which the centre diagram is also discrete, so that we may consider the cospan to be of the form 
$$\xymatrix@+1pc{\{A_i\}_{(i\in I)}\ar[r]^-{\phi:I\to J}&
\{B_{j}\}_{(j\in J)}&\{C_k\}_{(k\in K)}\ar[l]_-{\psi:J\leftarrow K}}.$$
Then colimit applied to this cospan is
$$\xymatrix@+1pc{\Sigma_{i\in I}A_i\ar[r]^-{colim(\phi)}&
\Sigma_{j\in J}B_j&\Sigma_{k\in K}A_k\ar[l]_-{colim(\psi)}},$$
where $colim(\phi)\cdot inj_i=inj_{\phi(i)}\ (i\in I)$ and $colim(\psi)\cdot inj_k=inj_{\psi(k)}\ (k\in K)$
 \begin{remark}
$\cspge$ is the result of freely adding wscc category structure to the graph $|\e|$  (a special case of this result was proved in \cite{RSWCT04}).
This means that diagrams in $|\e|$ may be written as expressions in the wscc structure of $\cspge$ with constants being the cospans of the form $\vec f$ for arrows $f$ of $|\e|$.
\end{remark}
\noindent Then colim preserves wscc expressions, so the colimit of any diagram may be written as an expression in $\cospe$.
 This is the compositionality of the calculation of colimits, mentioned in the title.

\subsection{The example of coequalizers}Consider the following cospan of diagrams in $\e$:
$$\xymatrix{\{A\}\ar[r]&
\{A\ar@<1ex>@{->}[r]^-{f}\ar@<-1ex>@{->}[r]_-{g}&B\}&\{B\}\ar[l]}.$$
The  cospan of diagrams may be pictured, as described above, as


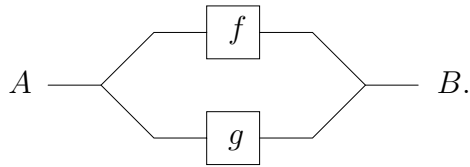
\centerline{{\tt\setlength{\unitlength}{1.00pt}
\begin{picture}(155,80)
\thinlines
\put(-15,37){$A$}
\put(147,37){$B.$}
\put(68,18){$g$}
\put(68,58){$f$}
\path(0,40)(20,40)(40,60)(60,60)(60,70)(80,70)(80,50)(60,50)(60,60)
\path(20,40)(40,20)(60,20)(60,30)(80,30)(80,10)(60,10)(60,20)
\path(80,20)(100,20)(120,40)(140,40)
\path(80,60)(100,60)(120,40)
\end{picture}}}
\noindent It is clear from the picture that the cospan may be expressed as the following composite in  $\cspge:$
$$\xymatrix@+1pc{
\{A\}\ \ar@{<->}[r]^-{\Delta}& \{A\}+\{A\} \ar@{<->}[r]^-{\{ f\}+\{g\}}&\{B\}+\{B\}\ar@{<->}[r]^-{\nabla}&\{B\}.\\ }$$
Applying colimit we see that the coequalizer of $f$ and $g$ may be expressed as the composite in $cospan(\e)$
$$\xymatrix{
A \ar@{<->}[r]^-{\Delta}& A+A\ar@{<->}[r]^-{\vec f+\vec g}&B+B\ar@{<->}[r]^-{\nabla}&B.\\ }$$
The composite of these three spans is the pushout Q of the following diagram in $E$
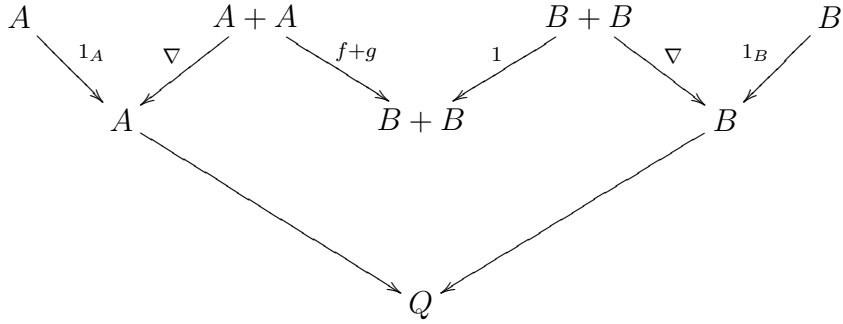
$$\xymatrix{
A \ar[rd]^-{1_{A}}&& A+A\ar[ld]_-{\nabla}\ar[rd]^-{f+g}&&B+B\ar[dl]_-{1}\ar[dr]^-{\nabla}&&B\ar[dl]_{1_B}\\
&A\ar[drdr]_-{}&&B+B&&B\ar[dldl]\\
&&&&&\\
&&&Q&&
}$$

\noindent It is not difficult to verify directly that $Q$ so defined is the coequalizer of $f$ and $g$.

\begin{example} By the same kind of reasoning the colimit of the diagram in \ref{compeqn} may be given by the expression of cospans in $\e$
$$(\epsilon_B+C)\cdot(\vec h+B+\vec k)\cdot (\Delta+B)\cdot\Delta\cdot(\nabla+\epsilon_A)\cdot(\vec f+B+\vec g+A)\cdot(A+\eta_B+A)\cdot \Delta.$$\end{example}
\begin{remark} If $\e$ has finite limits the functor $lim:(Graph/|\e|)^{op}\to\e$ extends to a functor
      $$ lim : \cspge \to  span(\e)$$
  which preserves the wscc structure.  This permits the compositional calculation of finite \emph{limits} in $\e$. In fact the equalizer of two arrows $\xymatrix{
A\ar@<1ex>@{->}[r]^-{f}\ar@<-1ex>@{->}[r]_-{g}&B}$ may be calculated by the same expression as that of the coequalizer above, but evaluated in $span(\e)$ rather than $cospan(\e)$.

\end{remark}

\subsection{Sketch of Proof of Theorem}
The main point to check in showing that colim is a monoidal functor is (a special case of) the following:

Consider a diagram $D$ of diagrams in $\e$ parametrized by a graph $G$;
 that is, a graph morphism $D : G \to   Graph/|\e|$. We can do two things.

\noindent (1) Calculate first the colimit of $D$ in $Graph/|\e|$ to obtain a diagram in $\e$ of which we may then take the colimit in $\e$,
       that is calculate $$colim_\e(colim_{Graph/|\e|}(D)).$$
(2) Calculate the colimit of 
$$\xymatrix{G\ar[r]^-{D}& Graph/\e\ar[r]^-{colim}&\e }$$
that is, calculate $colim_\e(colim_\e \cdot D)$.

\subsubsection{Lemma}
$colim_\e(colim_{Graph/|\e|}(D))\cong colim_\e(colim_\e \cdot D).$

\noindent Sketch of proof.

\noindent It suffices to show for any $X\in \e$ a bijection between cocones$$   colim_{Graph/|\e|}(D)\longrightarrow           X$$ and cocones
          $$ colim_{\e }\cdot D \longrightarrow X$$

But it is not hard to show that both of these are equivalent to a \textquotedblleft compatible family\textquotedblright\  of cocones
 $$D(g)\longrightarrow X\ \ \ \ (g \in G).$$

\subsubsection{A very special case}
Consider a diagram $D$ of diagrams in $\e$, namely $$D=\{\{A\} \ \{B\ C\}\}.$$
Then $$colim_\e(colim_{Graph/|\e|}(D))=colim_\e(\{A\ B\ C\})=A+B+C,$$
whereas $$colim_\e(colim_\e \cdot D)=colim_\e(\{A\ B+C\})=A+(B+C).$$
The lemma says exactly that the triple sum may be formed by repeated double sums, which has as a consequence the associative law for sums. It is clear that the general form of the lemma implies many further \textquotedblleft associative laws\textquotedblright\ - any two wscc expressions which yield the same diagram evaluate to the same in $cospan(\e)$.
 
\subsection{Example of Theorem}
A general cospan in $\cspge$ from $\varnothing$ to $\varnothing$ with centre the $D$ with may be constructed by taking the disjoint sum of all the arrows, and then equating vertices appropriately. This yields a formula for the general colimit of a finite diagram as follows. Let $\Sigma\scriptstyle{dom}$ denote the graph  $\sum_{\alpha\in D}\{dom(\alpha)\}$ and $\Sigma\scriptstyle{cod}$ denote the graph $\sum_{\alpha\in D}\{codom(\alpha)\}$. Let $\Sigma\alpha$ denote the graph $\sum_{\alpha\in D}\{\alpha\}$. Let $\Sigma\scriptstyle{ obj} $ denote the diagam consisting of all the objects in the $D$. Finally, let $i_{dom}$ and $i_{cod}$ denote the discrete cospans corresponding to the domain and codomain functions on the arrows of the graph parametrizing $D$. Then the cospan may be written
$$\xymatrix@+4pc{
\{\}\ar@{<->}[r]^-{\eta}&
\Sigma{\scriptstyle{dom}}+\Sigma{\scriptstyle{dom}}\ar@{<->}[r]^-{i_{cod}\cdot(\sum\alpha)+i_{dom}} &
\Sigma{\scriptstyle obj}+\Sigma{\scriptstyle obj}&
\{\}\ar@{<->}[l]_-{\epsilon}}.$$ Evaluating this formula instead in $\cospe $ gives the classical formula for colimits in terms of the coequalizer of two arrows from $\sum_{\alpha\in D}dom(\alpha)$ to $\sum_{\ A\in D}A$ ($\alpha$ arrow in $D$, $A$ object in $D)$.

\section{Limits and colimits of monoidal diagrams}
Systems in computer science are not usually  constructed from parts with "one input" and "one output", like arrows in a graph. Components have multiple inputs and outputs; that is, they are arrows in a monoidal graph.

\begin{definition}A \emph{monoidal graph} $(A,V,d_0,d_1)$ consists of a set $V$ of vertices, and a set $A$ of arcs and two functions $ d_0,d_1:A\longrightarrow\       V^*$ ($V^*$ the free monoid on $V$). A \emph{morphism of monoidal graphs} $\phi=(\phi_0,\phi_1)$  from $(A,Y,d_0,d_1)$ to $(B,W,d_0,d_1)$ consists of two functions $\phi_1:A\to B$ and $\phi_0:V\to W$ such that $\phi_0 d_0=d_0\phi_1,\phi_0d_1=d_1\phi_1$. We denote the category (actually a presheaf category) of monoidal graphs as $MonGraph$. There is an obvious notion then of a \emph{monoidal diagram} in a monoidal category since any monoidal category has an underlying monoidal graph.
\end{definition}

\begin{definition}Let $\e$ be a category with finite colimits, regarded as a monoidal category with sum as tensor. A cocone of a monoidal diagram $D$ to an object $X$ is a family of arrows $(q_{i}:A_{i}\longrightarrow X)\quad(A_{i}\text{ objects of the diagram $D$)}$
 such that for any arrow 
     $f:A_{i_1}+A_{i_2}+\cdots+A_{i_m}\to A_{j_1}+A_{j_1}+\cdots + … A_{j_n}$ in the diagram
     $$(q_{j_1}| q_{j_2}| q_{j_n}|\cdots  |q_{i_1})\cdot f = (q_{i_1}| q_{i_2}| q_{i_3}\cdots |q_{i_m}).$$ 
A colimit of  monoidal diagram $D$ is an object $C$ with a cocone $q$ from $D$ which is univeral; that is, any cocone to an object $X$ factors uniquely through $q$.
\end{definition}

\subsection{$\cspmge$}
Let $\e$ be a category with finite colimits, regarded as a monoidal category with sum as tensor. Then  $\cspmge$ is the full subcategory of $cospan(MonGraph/|\e|)$ whose objects are \emph{discrete} diagrams in $\e$. Just as with $\cspge$ we may picture arrows in $\cspmge$, the only difference being that components may have several input and output wires. Monoidal colimits are also calculable compositionally, in the algebra \( cospan(\mathbf{E})\), by a result analogous to \ref{theorem}. We look at one example only.
\subsection{Example}Consider the following cospan of monoidal diagrams $D$  in $\e$: the centre is the diagram with three objects $A$, $B$, $C$, and one arrow $f:A+C\to B+C$; the left hand side is the diagram $\{A\}$; the right hand side is the diagram $\{B\}$. Pictured the cospan is

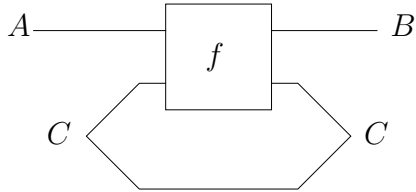
\centerline{{\tt\setlength{\unitlength}{1.0pt}
\begin{picture}(155,80)
\thinlines
\put(0,58){$A$}
\put(145,58){$B$}
\put(15,17){$C$}
\put(135,17){$C$}
\put(75,47){$f$}
\path(10,60)(60,60)(60,70)(100,70)(100,60)(140,60)
\path(60,60)(60,30)(100,30)(100,60)
\path(30,20)(50,40)(60,40)
\path(30,20)(50,0)(110,0)(130,20)(110,40)(100,40)
\end{picture}}}
\bigskip\noindent
From the picture we see immediately that this cospan of diagrams is expressible as a composite in $\cspmge$, namely
$$\xymatrix @+2pc {
\{A\}\ar@{<->}[r]^-{\{A\}+\eta_{\{C\}}\ }& \{A\}+\{C\}+\{C\}\ar@{<->}[r]^-{\{f\}+\{C\}}&\{B\}+\{C\}+\{C\}\ar@{<->}[r]^-{\{B\}+\epsilon_{\{C\}}}&\{B\}\\ 
}.$$

\noindent Applying monoidal colimit yields the fact that the monoidal colimit of the original diagram is the following composite in $cospan(\e)$:
$$\xymatrix{
A \ar@{<->}[r]^-{A+\eta_C }& A+C+C\ar@{<->}[r]^-{\vec f+C}&B+C+C\ar@{<->}[r]^-{B+\epsilon_{C}}&B\\ }.$$
Hence the colimit of the original diagram can be calculated as the  pushout below.
$$\xymatrix{
A \ar[rd]^-{inj}&& A+C+C\ar[ld]_-{A+\nabla}\ar[rd]^-{f+C}&&B+C+Cè\ar[dl]_-{B+C+C}\ar[dr]^-{B+\nabla }&&B\ar[dl]_{inj}\\
&A+C\ar[drdr]_-{}&&B+C+C&&B+C\ar[dldl]\\
&&&&&\\
&&&colim&&
}$$

\noindent
\medskip The colimit consists of orbits of $A+B+C$ under $f$. The pullback of the resulting cospan is the partial function obtained by iterating $f$.

\section{The Kleene Theorem}

\begin{theorem}\emph{(Kleene)} 

The languages recognized by finite state automata are the closure of singletons under union, concatenation and iteration.
\end{theorem} 

To prove this classical theorem the category \(\mathbf{E}\) we consider is \(\wp (\Sigma^*)\text{-}Cat\), categories enriched in \(\Sigma\text{-languages.}\) 
There is a composite of wscc functors$$\xymatrix{
Cspn(Graph/ \Sigma)\ar@{->}[r]^-{\Phi_1 }& Cspn(Graph/\wp(\Sigma^{*}))\ar@{->}[r]^-{\Phi_2}&               cospan(\wp(\Sigma^{*})\text{-}Cat). }$$

\noindent which takes a labelled graph (with input and output states) to the \(\wp(\Sigma\)*)-category whose homs are the languages traced out from the domain to the codomain. The existence of wscc functor $\Phi_1$ is implied by \cite{RSWCT04}, and of $\Phi_2$ (=colim) by \ref{theorem}.

This is already a Kleene-type theorem since, conceptually, the Kleene Theorem says that behaviour is an operation-preserving morphism from an algebra of systems to an algebra of possible behaviours, which implies that the perceived behaviours are the smallest class of possible behaviours closed under operations. In this case, the algebra of systems, that is the left-hand side, is generated as a wscc category by single labelled edges, and hence the image on the right-hand side is also generated as a wscc category by singleton languages.
 
However it is not the classical Kleene theorem, since the right-hand side does not consist of single languages and the wscc operations of $cospan(\wp(\Sigma^{*})\text{-}Cat)$ are not the Kleene operations. Further the functor does not lose internal states.

To obtain a theorem closer to the classical Kleene theorem we consider \emph{corelations}  between \(\wp(\Sigma\)*) categories, by which we mean\  cospans which are jointly bijective on objects. Then we compose the above wscc functor $\Phi_2 \Phi_1$ with a further wscc functor

             \[\Phi_3:cospan(\wp(\Sigma^*)\text{-} Cat)\longrightarrow corel(\wp(\Sigma^*)\text{-}Cat)\]

\noindent which uses the bijective-on-objects fully-faithful factorization to obtain from a cospan of \(\wp(\Sigma\)*)-categories a corelation of \(\wp(\Sigma\)*)-categories.

The final composite
             \[ \Phi_3\Phi_2\Phi_{1} :Cspn(Graph/\Sigma)\longrightarrow corel(\wp(\Sigma^*)\text{-}Cat)\]
\noindent takes a labelled graph with initial and final states to the category with objects only the initial and final states,  and homs the languages traced out.

To finish a proof of the classical Kleene theorem we need to show that the wscc operations in  $corel(\wp(\Sigma^*)\text{-}Cat)$, at the level of languages (homs), may be expressed in terms of union, concatenation and $()^*$.

Clearly the operation of tensor of corelations does not change the languages which occur as homs.
The problem is the composition.
But in a wscc category the  composition of two arrows may be expressed in terms of the tensor  and composition with the constants of the compact closed structure;
pictured, this is the fact that:

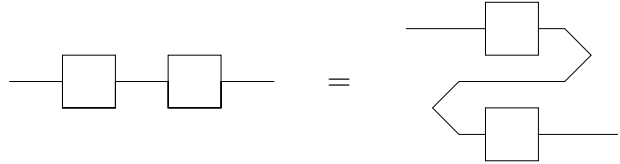
\centerline{{\tt\setlength{\unitlength}{1.0pt}
\begin{picture}(200,100)
\thinlines
\path(0,50)(20,50)(20,60)(40,60)(40,40)(20,40)(20,50)(20,40)(40,40)(40,50)(60,50)
\path(60,50)(60,60)(80,60)(80,40)(60,40)(60,50)(60,40)(80,40)(80,50)(100,50)
\put(120,47){$=$}
\path(150,70)(180,70)(180,80)(200,80)(200,60)(180,60)(180,70)
\path(200,70)(210,70)(220,60)(210,50)(170,50)(160,40)(170,30)(180,30)(180,40)(200,40)(200,20)(180,20)(180,30)
\path(200,30)(230,30)
\end{picture}}}

 So the general operation of composition in $corel(\wp(\Sigma^*)\text{-}Cat)$ may be reduced  to  the very special case of the colimit identifying two objects in a category.
So consider a $\wp(\Sigma^*)$-category $X$ containing two objects $x$ and $y$. The colimit category $X'=X/(x=y)$ has hom $X'(z,w)$ equal to $$X(z,w)\cup(X(z,x)\cup X (z,y))\cdot(X(x,x)\cup X(x,y)\cup X(y,x)\cup X(y,y))^*\cdot (X(x,w)\cup X (y,w)),$$ expressible using only Kleene operations. Hence the result.

This proof is very close to one of  the usual proofs of Kleene (if you strip the superstructure). Notice that the passage $\Phi_1$ permits the introduction of $\varepsilon$ moves; that is, homs which consist of the empty word.  The superstructure has the advantage of suggesting needed generalizations, for example, to parallelism.

\section{Comments}

The theorem we have described concerns calculating colimits as objects, not as functors. We have not shown the  compositionality of morphisms between colimits. We believe that this is connected with the algebra of span and cospan as symmetric monoidal bicategories, rather than as categories. We have made an initial progress in understanding this question in \cite{MSW07}, by considering a very special case, where we identify the role of "2-separable object".



\end{document}